\documentclass[11pt]{article}
\usepackage{epsfig,amssymb}
\begin{document}
\title{Asymptotic behavior of discrete holomorphic maps $z^c$,  $log(z)$ and discrete Painleve transcedents.}
\author{{\Large Agafonov S.I. } \\
\\
Institut f\"{u}r Algebra und Geometrie\\
Martin-Luther-Universit\"{a}t
Halle-Wittenberg \\ D-06099 Halle (Saale), Germany\\
e-mail: {\tt agafonov@mathematik.uni-halle.de} }
\date{}

\newtheorem{theorem}{Theorem}
\newtheorem{proposition}{Proposition}
\newtheorem{lemma}{Lemma}
\newtheorem{corollary}{Corollary}
\newtheorem{definition}{Definition}
\newtheorem{conjecture}{Conjecture}

\pagestyle{plain}

\maketitle

\begin{abstract}
It is shown that discrete analogs of $z^{c }$ and ${\rm log}(z)$
have the same asymptotic behavior as their smooth counterparts.
These discrete maps are described in terms of special solutions of
discrete Painleve-II equations, asymptotics of these solutions
providing the behaviour of discrete $z^{c }$ and ${\rm log}(z)$
at infinity.
\end{abstract}

\section{Introduction}

In this paper we discuss a very interesting and rich object:
circle patterns mimicking the holomorphic maps $z^{c }$ and ${\rm
log}(z)$. Discrete $z^{2 }$ and ${\rm log}(z)$ were guessed by
Schramm and Kenyon (see \cite{www}) as examples from the class of
circle patterns with the combinatorics of the square grid
introduced by Schramm in \cite{Schramm}. This class opens a new
page in the classical theory of circle packings, enjoying new
golden age after Thurston's idea \cite{T} about approximating the
Riemann mapping by circle packings: it turned out that these
circle patterns are governed by discrete integrable equation (the
stationary Hirota equation), thus providing one with the whole
machinery of the integrable system theory (\cite{BPD}).

The striking analogy between circle patterns and holomorphic maps
resulted in the development of discrete analytic function theory
(for a good survey see\cite{DS}). Classical circle packings
comprised of disjoint open disks were later generalized to circle
patterns where the disks may overlap. Discrete versions of
uniformization theorem, maximum principle, Schwarz's lemma and
rigidity properties and  Dirichlet principle were established
(\cite{MR},\cite{H},\cite{Schramm}).

Different underlying combinatorics were considered: Schramm
introduced square grid circle patterns, generalized by Bobenko
and Hoffmann to hexagonal patterns with constant intersection
angles in \cite{BH}, hexagonal circle patterns with constant
multi-ratios were studied by Bobenko, Hoffman and Suris in
\cite{BHS}.

The difficult question of convergence was settled by Rodin and
Sullivan \cite{RS} for general circle packings, He and Schramm
\cite{HS} showed that the convergence is $C^{\infty}$ for
hexagonal packings, the uniform convergence for square grid
circle patterns  was established by Schramm \cite{Schramm}.

On the other hand not very many examples are known: for circle
packings with the hexagonal combinatorics the only explicitly
described examples are Doyle spirals, which are discrete
analogues of exponential maps \cite{Doy}, and conformally
symmetric packings, which are analogues of a quotient of Airy
functions \cite{BHConf}. For patterns with overlapping circles
more examples are constructed: discrete versions of ${\rm exp}
(z)$, ${\rm erf}(z)$ (\cite{Schramm}),  $z^c$, $ {\rm log} (z)$
(\cite{AB}) are constructed for patterns with underlying
combinatorics of the square grid; $z^c$, ${\rm log}(z)$ are also
described for hexagonal patterns with both multi-ratio
(\cite{BHS}) and constant angle (\cite{BH}) properties.

Discrete $z^c$ is not only a very interesting example in discrete
conformal geometry. It has mysterious relationships to other
fields. It is constructed via some discrete isomonodromic problem
and is governed by discrete Painlev\'e II equation
(\cite{AB},\cite{N}), thus giving geometrical
 interpretation thereof. Its linearization defines Green's function on critical
graphs (see \cite{BMS}) found in \cite{K} in the frames of the
theory of Dirac operator. Moreover, it seems to be a rather
important tool for investigation of more general circle patterns
and discrete minimal surfaces (see \cite{B} for a brief survey and
\cite{BHSp} for more details).

\subsection{Circle patterns and discrete conformal maps}

\begin{figure}[th]
\begin{center}
\epsfig{file=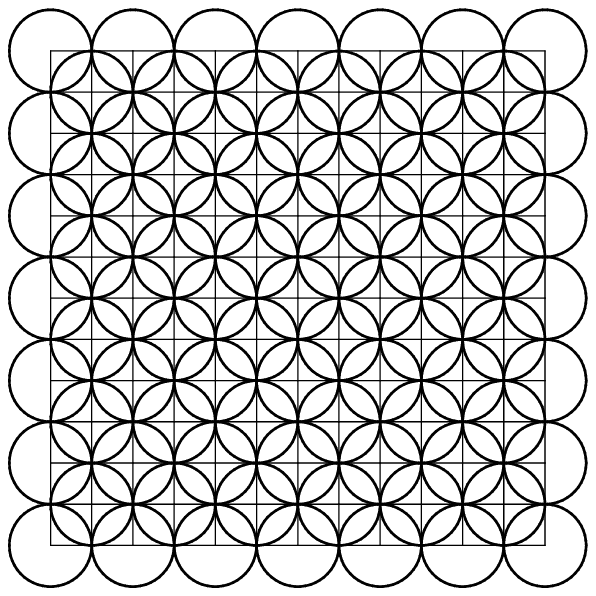,width=45mm}
\begin{picture}(20,50)
\put(3,52){\vector(1,0){20}} \put(7,58){ \it  \huge  f}
\end{picture}
\epsfig{file=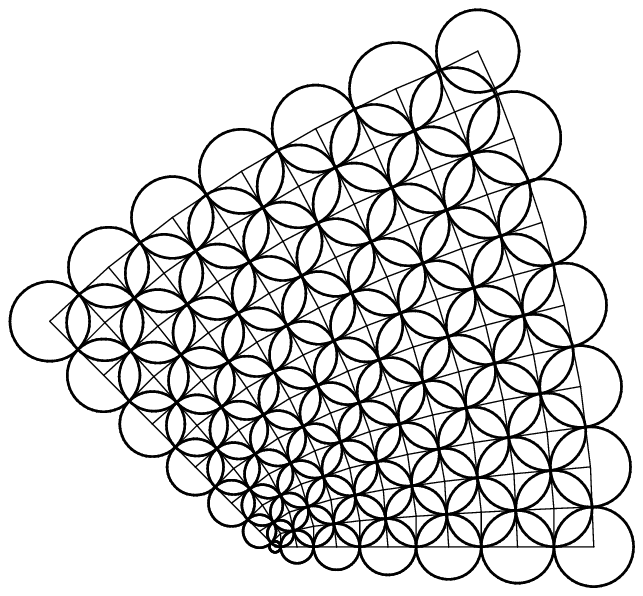,width=50mm} \caption{Schramm's circle
patterns as discrete conformal map.  Shown is the discrete
version of the holomorphic mapping $z^{3/2}$.} \label{CPMap}
\end{center}
\end{figure}

To visualize the analogy between Schramm's circle patterns and
conformal maps, consider regular patterns composed of unit
circles  and suppose that the radii are being deformed so as to
preserve the orthogonality of neighboring circles and the
tangency of half-neighboring ones. Discrete maps taking
 intersection points of
the unit circles of the standard regular patterns to the
respective points of the deformed patterns mimic classical
holomorphic functions, the deformed radii being analogous to
$|f'(z)|$ (see Fig.~\ref{CPMap}).

It is easy to show that the lattice comprised of the centers of
circles of  Schramm's pattern and their intersection points is a
special discrete conformal mapping (see Definition \ref{DCM}
below). The latter were introduced in \cite{BPdis} in the frames
of discrete integrable geometry, originally without any relation
to circle patterns.

\begin{definition} A map  $f\ : \ {\mathbb Z^2\ \rightarrow \ {\mathbb R^2}={\mathbb C}}$
is called a discrete conformal map if all its elementary
quadrilaterals are conformal squares, i.e., their cross-ratios
are equal to -1:

$$
q(f_{n,m},f_{n+1,m},f_{n+1,m+1},f_{n,m+1}):=
$$
\begin{equation}
\frac{(f_{n,m}-f_{n+1,m})(f_{n+1,m+1}-f_{n,m+1})}
{(f_{n+1,m}-f_{n+1,m+1})(f_{n,m+1}-f_{n,m})}=-1. \label{q}
\end{equation}
\label{DCM}
\end{definition}
This
definition is motivated by the following properties:\\
1) it is
M\"obius invariant,\\
2) a smooth map $f:\ D\ \subset {\mathbb C} \to {\mathbb C}$ is
conformal (holomorphic or antiholomorphic) if and only if
$$
\lim _{\epsilon \to 0}q(f(x,y),f(x+\epsilon ,y)f(x+\epsilon
,y+\epsilon )f(x,y+\epsilon ))=-1
$$
for all $(x,y)\in D$.

\begin{figure}[ht]
\begin{center}
\epsfig{file=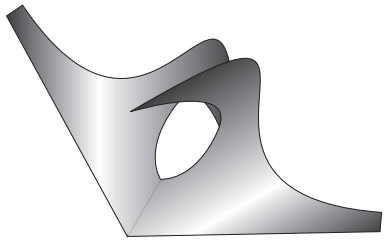, width=50mm}
\end{center}
\caption{Surface glued of quadrilaterals of immersed but
non-embedded discrete map. } \label{surface}
\end{figure}

\begin{definition}
A discrete conformal map $f_{n,m}$ is called embedded if
interiors of different elementary quadrilaterals
$(f_{n,m},f_{n+1,m},f_{n+1,m+1},f_{n,m+1})$ do not intersect.
\end{definition}
\noindent The condition for discrete conformal map to be embedded
can be relaxed as follows.
\begin{definition}
A discrete conformal map $f_{n,m}$ is called an immersion if
interiors of adjacent elementary quadrilaterals
$(f_{n,m},f_{n+1,m},f_{n+1,m+1},f_{n,m+1})$ are disjoint.
\end{definition}
To illustrate the difference between the immersed and embedded
discrete conformal maps, let us imagine that the elementary
quadrilaterals of the map are made of elastic inextensible
material and glued along the corresponding edges to produce
 a surface with a border. If this surface is immersed it is locally flat. Being dropped down it will not have folds.
At first sight it seems to be sufficient to guarantee
embeddedness, provided $Z_{n,0}\to \infty$ and $Z_{0,m}\to \infty$
as $n \to \infty$ (which follows from $ Z^{c }_{n,0}=k(c)n^{c
}\left(1+O\left(\frac{1}{n^2} \right) \right), \ \ n \to \infty $
 (\cite{AB}). But a surface with such properties still may have
some limit curve with self-intersections thus giving overlapping
quadrilaterals. Hypothetical example of such a surface is shown
in Fig.~\ref{surface}.

Equation (\ref{q}) can be supplemented with the following
nonautonomous constraint:
\begin{equation}
c f_{n,m}=2n\frac{(f_{n+1,m}-f_{n,m})(f_{n,m}-f_{n-1,m})}
{(f_{n+1,m}-f_{n-1,m})}+2m\frac{(f_{n,m+1}-f_{n,m})(f_{n,m}-f_{n,m-1})}
{(f_{n,m+1}-f_{n,m-1})}. \label{c}
\end{equation}
This constraint, as well as its compatibility with (\ref{q}), is
derived from some monodromy problem (see \cite{AB} for detail).
Let us assume $0<c < 2$ and denote ${\mathbb Z^2_{+}}=\{ (n,m) \in
{\mathbb Z^2}: n,m \ge 0 \}.$ Motivated by the asymptotics of the
constraint (\ref{c}) at $n,m \rightarrow \infty$ and the
properties

$$
z^c({\mathbb R_+})  \in {\mathbb R_+}, \ \ z^c(i{\mathbb R_+}) \in
e^{c \pi i/2}{\mathbb R_+}
$$
of the holomorphic mapping $z^c$ we use the following definition
(\cite{BPD}) of the "discrete" $z^c$.

\begin{definition}
The discrete conformal map $Z^c\ : \ {\mathbb Z^2_+\ \rightarrow \
{\mathbb C}}, \  0<c < 2\ $ is the solution of (\ref{q}),(\ref{c})
with the initial conditions

\begin{equation}
Z^c(0,0)=0, \ \ Z^c(1,0)=1, \ \ Z^c(0,1)=e^{c \pi i/2}.
\label{initial}
\end{equation}
\label{def}
\end{definition}
Obviously, $Z^c(n,0)\in {\mathbb R_{+}}$ and
 $Z^c(0,m)\in
e^{c \pi i/2} ({\mathbb R_{+}})$  for any  $  n,m
\in {\mathbb N}$ .\\
A limit of rescaled $Z^c$ as $c$ approaches 2 gives discrete
$Z^2$. Discrete ${\rm Log}$
  is defined as dual to $Z^2$ (see subsection \ref{s.Z2} for detail).
\begin{theorem}(\cite{A}) \label{emb}
The discrete maps ${\rm Log}$ and $Z^{c }$ for $0<c \le 2$ are
embedded.
\end{theorem}

Given initial data $f_{0,0}=0$, $f_{1,0}=1$,  $f_{0,1}=e^{i \alpha
}$ with $\alpha \in \mathbb R$,  constraint (\ref{c}) allows one
to compute $f_{n,0}$ and $f_{0,m}$ for all $n,m \ge 1.$ Now using
equation (\ref{q}) one can successively compute $f_{n,m}$ for any
$n,m \in {\mathbb N}$. It turned out that all edges at the vertex
$f_{n,m}$ with $n+m=0 \ ({\rm mod} \ 2)$ are of the same length
\begin{equation}\label{equi}
|f_{n+1,m}-f_{n,m}|=|f_{n,m+1}-f_{n,m}|=|f_{n-1,m}-f_{n,m}|=|f_{n,m-1}-f_{n,m}|
\end{equation}
and  all angles between the neighboring edges at the vertex
$f_{n,m}$ with $n+m=1 \ ({\rm mod} \ 2)$  are equal to $\pi /2.$
Thus for any $n,m: \ n+m=0 \ ({\rm mod }\ 2)$ the points
$f_{n+1,m},$ $f_{n,m+1},$ $f_{n-1,m},$ $f_{n,m-1}$ lie on the
circle with the center $f_{n,m}$. All such circles form a circle
pattern of Schramm type (see \cite{Schramm}), i.e. the circles of
neighboring quadrilaterals intersect orthogonally and the circles
of half-neighboring quadrilaterals with common vertex are
tangent. Consider the sublattice $\{n,m: \ n+m=0 \ ({\rm mod} \
2)\}$ and denote by $\mathbb V$ its quadrant
\begin{equation}\label{V}
{\mathbb V}=\{z=N+iM:\ N,M \in {\mathbb Z^2}, M \ge |N| \},
\end{equation}
where
$$
N=(n-m)/2, \ \ M=(n+m)/2.
$$
We will use complex labels $z=N+iM$ for this sublattice. Denote by
$C(z)$ the circle of the radius
\begin{equation}
R(z)=|f_{n,m}-f_{n+1,m}|=|f_{n,m}-f_{n,m+1}|=|f_{n,m}-f_{n-1,m}|=|f_{n,m}-f_{n,m-1}|
\label{rmap}
\end{equation}
with the center at $f_{N+M,M-N}=f_{n,m}.$ Then $Z^c$ is immersed
if and only if $R(z)$ satisfy the following equation
\begin{equation}
R(z)^2=\frac{\left(
\frac{1}{R(z+1)}+\frac{1}{R(z+i)}+\frac{1}{R(z-1)}+\frac{1}{R(z-i)}
\right) R(z+1)R(z+i)R(z-1)R(z-i) }{ R(z+1)+R(z+i)+R(z-1)+R(z-i)  }
\label{lnR}
\end{equation}
The proof of Theorem \ref{emb}
 is based on the analysis of the solution to the
following equations:
\begin{equation}
\begin{array}{l}
R(z)R(z+1)(-2M-c
)+R(z+1)R(z+1+i)(2(N+1)-c )+\\
\qquad R(z+1+i)R(z+i)(2(M+1)-c )+R(z+i)R(z)(-2N-c )=0,
\end{array}
\label{square}
\end{equation}
\begin{equation}
\begin{array}{l}
(N+M)(R(z)^2-R(z+1)R(z-i))(R(z+i)+R(z+1))+\\
\qquad (M-N)(R(z)^2-R(z+i)R(z+1))(R(z+1)+R(z-i))=0,
\end{array}
\label{Ri}
\end{equation}
$$
(N+M)(R(z)^2-R(z+i)R(z-1))(R(z-1)+R(z-i))+
$$
\begin{equation}
(M-N)(R(z)^2-R(z-1)R(z-i))(R(z+i)+R(z-1))=0, \label{Le}
\end{equation}
$$
(N+M)(R(z)^2-R(z+i)R(z-1))(R(z+1)+R(z+i))+
$$
\begin{equation}
(N-M)(R(z)^2-R(z+1)R(z+i))(R(z+i)+R(z-1))=0, \label{Up}
\end{equation}
which are compatible with (\ref{lnR}). Namely, $Z^c$ is immersed
iff the solutions of these equations in ${\mathbb V}$ with
\begin{equation}
R(0)=1, \ \ \ R(i)=\tan \frac{c \pi}{4}. \label{Rinitial}
\end{equation}  is
positive (see \cite{AB}). More delicate global property of
embeddedness follows from the following proposition.
\begin{proposition}(\cite{A})
If for  a  solution R(z) of (\ref{square}, \ref{Ri}) with $c \ne
1$ and initial conditions (\ref{Rinitial})  holds
\begin{equation}
 R(z)>0, \ \ (c -1)(R(z)^2 - R(z-i)R(z+1))\ge 0
\label{sign}
\end{equation}
  in interior vertices of ${\mathbb V}$, then $Z^c$ is embedded. \label{convex}
\end{proposition}

\subsection{The discrete maps $Z^2$ and $\rm Log $. Duality}
\label{s.Z2}

Definition \ref{def} was given for $0<c <2.$ For $c <0$ or $c
>2$, the radius $R(1+i)=c /(2-c)$ of the corresponding
circle patterns (found as a solution to equations for $R(z)$)
 becomes negative and some elementary
quadrilaterals around $f_{0,0}$ intersect. But for $c =2$, one
can renormalize the initial values of $f$ so that the
corresponding map remains an immersion. Let us consider $Z^{c }$,
with $0<c <2$,
 and make the following renormalization
for the corresponding radii:  $R \to \frac{2-c }{c }R.$ Then as
$c \to 2-0$  we have
$$
R(0)=\frac{2-c }{c }\to +0, \ \  R(1+i)=1, \ \ R(i)=\frac{2-c }{c
}\tan \frac{c \pi}{4}\to \frac{2}{\pi}.
$$
\begin{definition}(\cite{AB})
$Z^{2}\ : \ {\mathbb Z^2_+\ \rightarrow \ {\mathbb R^2}={\mathbb
C}}$ is the solution of (\ref{q}), (\ref{c}) with $c =2$ and the
initial conditions
$$
Z^{2}(0,0)=Z^{2}(1,0)=Z^{2}(0,1)=0, \ \ Z^{2}(2,0)=1, \ \
Z^{2}(0,2)=-1, \ \ Z^{2}(1,1)=i\frac{2}{\pi}.
$$
\label{def2}
\end{definition}

\begin{figure}[ht]
\begin{center}
\epsfig{file=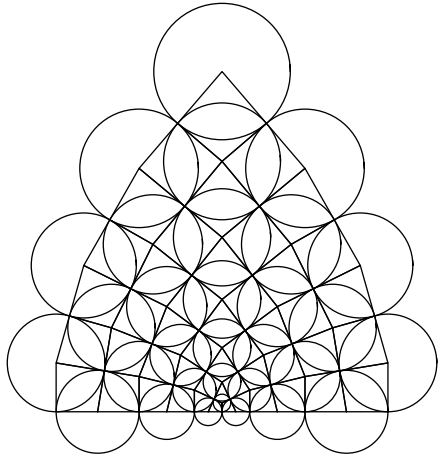,width=60mm} \caption{Discrete $Z^2$.}
\label{f.Z2}
\end{center}
\end{figure}

In this definition, equations (\ref{q}),(\ref{c}) are understood
to be regularized through multiplication by their denominators.
Note that for the radii on the border one has $R(N+iN)=N.$  \\
If  $R(z)$ is a solution to (\ref{lnR}) and therefore defines
some immersed circle patterns, then $\tilde R(z)=\frac{1}{R(z)}$
also solves (\ref{lnR}). This reflects the fact that for any
discrete conformal map $f$ there is {\it dual discrete conformal
map} $f^*$ defined by (see \cite{BPD})
\begin{equation}                        \label{duality}
f^*_{n+1,m}-f^*_{n,m}=-\frac{1}{{f_{n+1,m}}-{f_{n,m}}}, \ \
f^*_{n,m+1}-f^*_{n,m}=\frac{1}{{f_{n,m+1}}-{f_{n,m}}}.
\end{equation}
The smooth limit of the duality (\ref{duality}) is
$$
({f^*})'=-\frac{1}{f'}.
$$
The dual of $f(z)=z^2$ is, up to a constant, ${f^*(z)}=\log z.$
Motivated by this observation, we define  the discrete logarithm
as the discrete map dual to $Z^2$, i.e. the map corresponding to
the circle pattern with radii
$$
R_{\rm{Log}}(z)=\frac{1}{R_{Z^2}(z)},
$$
where $R_{Z^2}$ are the radii of the circles for $Z^2.$ Here one
has $R_{\rm{Log}}(0)=\infty$, i.e. the corresponding circle is a
straight line. The corresponding constraint  (\ref{c})
 can be also derived as a limit. Indeed, consider the map
 $g=\frac{2-c }cZ^{c }-\frac{2-c }c.$ This map satisfies (\ref{q})
 and the constraint
$$
c \left( g_{n,m}+ \frac{2-c }c\right)=
2n\frac{(g_{n+1,m}-g_{n,m})(g_{n,m}-g_{n-1,m})}
{(g_{n+1,m}-g_{n-1,m})}+
2m\frac{(g_{n,m+1}-g_{n,m})(g_{n,m}-g_{n,m-1})}
{(g_{n,m+1}-g_{n,m-1})}.
$$
Keeping in mind the limit procedure used do determine $Z^2$, it
is natural to define the discrete analogue of $\log z$ as the
limit of $g$ as $c \to +0$. The corresponding constraint becomes
\begin{equation}
1=n\frac{(g_{n+1,m}-g_{n,m})(g_{n,m}-g_{n-1,m})}
{(g_{n+1,m}-g_{n-1,m})}+m\frac{(g_{n,m+1}-g_{n,m})(g_{n,m}-g_{n,m-1})}
{(g_{n,m+1}-g_{n,m-1})}. \label{cln}
\end{equation}

\begin{figure}[ht]
\begin{center}
\epsfig{file=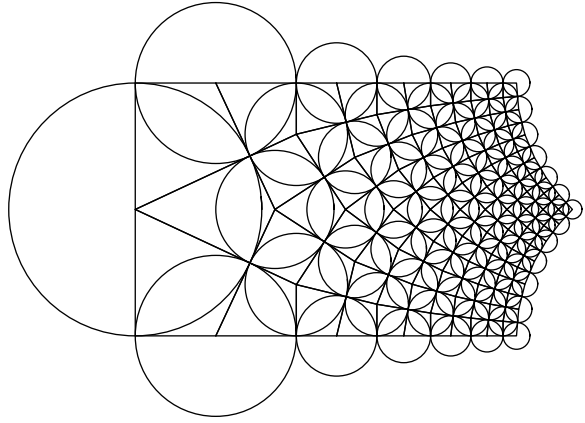,width=80mm} \caption{Discrete {\rm Log}.}
\label{f.Log}
\end{center}
\end{figure}
\begin{definition} (\cite{AB}) $\rm{Log}$  is the map
$\rm{Log} : \ {\mathbb Z^2_+ \rightarrow \ {\mathbb R^2}=\bar
{\mathbb C}}$ satisfying (\ref{q}) and (\ref{cln}) with the
initial conditions
$$
\rm{Log}(0,0)=\infty , \ \rm{Log}(1,0)=0, \ \rm{Log}(0,1)=i\pi ,
$$
$$
 \ \rm{Log}(2,0)=1, \ \rm{Log}(0,2)=1+i\pi ,
  \ \rm{Log}(1,1)=i\frac{\pi}{2}.
$$
\label{defLn}
\end{definition}
The circle patterns corresponding to the discrete conformal
mappings $Z^2$ and $\rm{Log}$ were conjectured
 by O.~Schramm and R.~Kenyon (see \cite{www}), in \cite{AB} it was proved that
they are immersed.

\subsection{ $Z^c$ and discrete Painlev\'e equations} The main tool
to establish all the above mentioned properties of $R(z)$ was a
special case of discrete Painlev\'e-II equation:
\begin{equation}
\label{dPII}
(n+1)(u_n^2-1)\left(\frac{u_{n+1}-iu_n}{i+u_nu_{n+1}}\right)-
 n(u_n^2+1)\left(\frac{u_{n-1}+iu_n}{i+u_{n-1}u_n}\right) =
c u_n.
\end{equation}
The embedded $Z^c$ corresponds to the unitary solution
$u_n=e^{i\alpha_n}$ of this equation with $u_0=e^{ic \pi/4}$
$0<\alpha _n < \pi/2.$ ($\alpha_n$ is define by
$f_{n,n+1}-f_{n,n}=e^{2i\alpha _n}(f_{n+1,n}-f_{n,n}).$)

In \cite{A} was proved:
$$
\lim _{n\to \infty}Z^c _{n,m}=\infty, \ \ \lim _{m\to \infty}Z^c
_{n,m}=\infty.
$$
In the present paper we use (\ref{Ri}) and (\ref{square}) (which
form of course a system of Painlev\'e type) to establish the more
accurate result
$$R(N_0+iM) \simeq K(c)M^{c-1} \ \ \ {\rm as} \ \ \ M\to
\infty
$$
and corresponding asymptotics for $Z^c$ and $u_n$. Thus asymptotic
behavior of discrete $z^c$ and $log(z)$ is exactly that of their
smooth counterparts.

 For smooth Painlev\'e equations similar asymptotic problems have been studied in the frames of the
isomonodromic deformation method \cite{IN}. In particular,
connection formulas were derived. These formulas describe the
asymptotics of solutions  for $n \to \infty$ as a function of
initial conditions. Some discrete Painlev\'e equations were
studied in that framework in \cite{FIK}). The geometric origin of
our equations permits us to find our asymptotics by bare-handed
approach, studying linearized equations.

\section{Asymptotics of discrete $z^c$ and ${\rm log}(z)$}

To treat $z^c$ and ${\rm log}(z)$ on equal footing we agree that
the case $c=0$ in equations corresponds to discrete ${\rm
log}(z)$.  For the edges of unit squares with the vertices in $
{\mathbb V}$ defined by (\ref{V}) we introduce $X$ and $Y$ via the
radius ratios:
\begin{equation} \label{t-def}
\frac{1+X_{N,M}}{1-X_{N,M}}=\frac{R_{N+1,M}}{R_{N,M}}, \ \ \ \
\frac{1+Y_{N,M}}{1-Y_{N,M}}=\frac{R_{N,M}}{R_{N,M-1}}.
\end{equation}
In these variables equations (\ref{Ri}) and (\ref{square}) read
as:
\begin{equation} \label{Ri-t}
(M-N)\frac{X_{N,M}+Y_{N,M+1}}{1-X_{N,M}Y_{N,M+1}}+(M+N)\frac{X_{N,M}-Y_{N,M}}{1+X_{N,M}Y_{N,M}}=0,
\end{equation}
\begin{equation} \label{square-t}
(M-N)\frac{Y_{N,M+1}-X_{N,M}}{1-X_{N,M}Y_{N,M+1}}+(M+N+1)\frac{X_{N,M+1}+Y_{N,M+1}}{1+X_{N,M+1}Y_{N,M+1}}=c-1.
\end{equation}
Moreover, $X$ and $Y$ satisfy
\begin{equation} \label{inR}
\frac{X_{N,M}+Y_{N,M+1}}{1-X_{N,M}Y_{N,M+1}}=\frac{X_{N-1,M}+Y_{N,M}}{1-X_{N-1,M}Y_{N,M}},
\end{equation}
\begin{equation} \label{outR}
\frac{X_{N,M}+Y_{N+1,M+1}}{1+X_{N,M}Y_{N+1,M+1}}=\frac{X_{N,M+1}+Y_{N,M+1}}{1+X_{N,M+1}Y_{N,M+1}},
\end{equation}
where (\ref{inR}) is equivalent to (\ref{lnR}) and (\ref{outR})
is the compatibility condition.

Conditions (\ref{sign}), which hold for discrete $z^c$ and ${\rm
log}(z)$, turn out to be so restrictive for the corresponding
solutions $X_{N,M},Y_{N,M}$ of (\ref{Ri-t}), (\ref{square-t}),
that they allow one to compute their asymptotic behavior. For
definiteness we consider the case $c>1$.
\begin{lemma} \label{estimate}
For the solution $X_{N,M},Y_{N,M}$ of
(\ref{Ri-t}),(\ref{square-t}) corresponding to $z^c$ in ${\mathbb
V}$ with $c>1$ holds true:
\begin{equation} \label{estimation}
-\frac{c-1}{M-N}\le X_{N,M} \le \frac{c-1}{M+N} \ \ \ {\rm and}\ \
0 \le Y_{N,M+1} \le \frac{c-1}{M+N}+\frac{2(c-1)}{M-N}.
\end{equation}
\end{lemma}
\noindent {\it Proof:} Note that for the studied solutions
$R\ge0$ and therefore $-1\le X_{N,M}\le 1$ and $-1\le
Y_{N,M}\le1$. Let us denote for brevity $X_{N,M}$ by $X$,
$X_{N,M+1}$ by $\bar X$, $Y_{N,M}$ by $Y$, $Y_{N,M+1}$ by $\bar
Y$, and $R_{N,M}$ by $R_1$, $R_{N+1,M}$ by $R_2$, $R_{N+1,M+1}$ by
$R_3$, $R_{N,M+1}$ by $R_4$, $R_{N,M-1}$ by $R_5$. Then
(\ref{sign}) together with (\ref{Ri}),(\ref{Le}),(\ref{Up}) imply
\begin{equation} \label{R-inequalities}
R_1^2\le R_2R_4, \ \ R_1^2\ge R_2R_5, \ \ R_3^2\ge R_2R_4, \ \
R_4^2\ge R_1R_3, \ \ R_2^2\le R_1R_3.
\end{equation}
Rewriting the first inequality as $\frac{R_1}{R_2}\le
\frac{R_4}{R_1}$ and taking into account $1+X\ge 0$ we have
\begin{equation}\label{t1+t2}
X+\bar Y \ge 0.
\end{equation}
Similarly the second  inequality of (\ref{R-inequalities}) infers
$Y-X\ge 0$ or after shifting
\begin{equation}\label{t2t6}
\bar Y \ge \bar X.
\end{equation}
Combining the first and the third inequalities of
(\ref{R-inequalities}) one gets $R_3\ge R_1$ or
\begin{equation}\label{t2+t6}
\bar X +\bar Y \ge 0,
\end{equation}
which is equivalent to
\begin{equation}\label{t5+t1}
X +\bar Y_{N+1,M+1} \ge 0
\end{equation}
due to (\ref{outR}). Similarly the forth and the fifth imply
$R_4\ge R_2$ or
\begin{equation}\label{t2t1}
\bar Y  \ge X.
\end{equation} Comparing (\ref{t2+t6}) with (\ref{t2t6}) one gets $\bar Y\ge 0$.
Rewriting (\ref{square}) as
$$
\frac{R_3}{R_1}=\frac{(2N+c)R_4+(2M+c)R_2}{(2(M+1)-c)R_4+(2(N+1)-c)R_2}
$$
and taking into account $\frac{R_3}{R_1}\ge1$ and $c \le 2 $ we
can estimate $\frac{R_2}{R_4}$ in ${\mathbb V}$:
$$
\frac{R_2}{R_4}\ge\frac{1-\frac{c-1}{M-N}}{1+\frac{c-1}{M-N}}
$$
which reads as
$$
\bar Y-X \le (1-\bar Y X)\epsilon _{\scriptscriptstyle N,M}, \ \ \
\epsilon _{\scriptscriptstyle N,M}=\frac{c-1}{M-N}.
$$
As $(1-\bar Y X)\le 2$ we have
\begin{equation}\label{t2-t1e}
\bar Y-X \le 2\epsilon _{\scriptscriptstyle N,M}.
\end{equation}
Similarly we can solve (\ref{square}) with respect to
$\frac{R_4}{R_2}$:
$$
\frac{R_4}{R_2}=\frac{(2M+c)R_1+(-2(N+1)+c)R_3}{(-2N-c)R_1+(2(M+1)-c)R_3}.
$$
Note that for $M > N$ holds $(-2N-c)R_1+(2(M+1)-c)R_3\ge 0$ as
$R_3\ge R_1$, and $c\le 2$. This together with $R_4\ge R_2$ and
$(1-\bar Y \bar X)\le 2$ gives after some calculations
\begin{equation}\label{t2+t6d}
\bar X+\bar Y \le 2\delta _{\scriptscriptstyle N,M},
\end{equation}
with $\delta _{\scriptscriptstyle N,M}=\frac{c-1}{M+N+1}$. Now
inequalities (\ref{t1+t2}),(\ref{t2t1}),(\ref{t2-t1e}) yield
$$
-\epsilon _{\scriptscriptstyle N,M}\le X
$$
and (\ref{t2t6}),(\ref{t2+t6}),(\ref{t2+t6d}) imply
$$
\bar X \le \delta _{\scriptscriptstyle N,M},
$$
which gives the first inequality of (\ref{estimation}) after
shifting backwards from $M+1$ to $M$ . Using
(\ref{t2-t1e}) we easily get the second one. $\Box$\\

The estimations obtained allow one to find asymptotic behavior
for the radius-function
in $M-$direction.
\begin{theorem}\label{M-ass}
For the solution $R_{N,M}$ of (\ref{Ri}),(\ref{square})
corresponding to discrete $Z^c$ and $\rm{Log}$ in ${\mathbb V}$
holds true:
\begin{equation}\label{as2}
R_{N_0,M} \simeq K(c)M^{c-1} \ \ \ {\rm as} \ \ \ M\to \infty,
\end{equation}
with constant $K(c)$ independent of $N_0$.
\end{theorem}
\noindent {\it Proof:} Because of the duality
$$R_{z^{2-c}}=\frac{1}{R_{z^c}}$$ it is enough to consider the case $c > 1$.
Let us introduce $n=M-N_0$, $x_n=X_{N_0,M}$, $y_n=Y_{N_0,M}$.
Then for large $n$ Lemma \ref{estimate} allows one to rewrite
equations (\ref{Ri-t}),(\ref{square-t}) as
\begin{equation} \label{liearize}
\begin{array}{l}
x_{n+1}=5x_n-2y_n+\frac{c-1}{n}+U_{N_0}(n),\\
y_{n+1}=-2x_n+y_n+V_{N_0}(n),
\end{array}
\end{equation}
where $U_{N_0}(n)$ and $V_{N_0}(n)$ are defined by discrete $Z^c$
and satisfy
\begin{equation} \label{restAB}
U_{N_0}(n) < C_1\frac{1}{n^2}, \ \ \ V_{N_0}(n) <
C_2\frac{1}{n^2},
\end{equation}
for natural $n$ and some constants $C_1,C_2$ (depending on $N_0$).
Thus the solution $(x_n,y_n)$, corresponding to discrete $Z^c$ is
a special solution of linear non-homogeneous system
(\ref{liearize}) having the order $\frac{1}{n}$ for large $n$.
The eigenvalues of the system matrix are positive numbers
$\lambda _1=3-2\sqrt 2<1$ and $\lambda _2=3+2\sqrt 2>1.$ In the
diagonal form system (\ref{liearize}) takes the form:
\begin{equation} \label{diagonal}
\varphi _{n+1}= A\varphi _n + s_n+r_n
\end{equation}
with
$$
A=\left (\begin{array}{cc} \lambda _1 & 0 \\
0& \lambda _2
\end{array} \right ), \ \ \varphi _n=\left (\begin{array}{cc} 1, & 1 \\
1+\sqrt 2 ,& 1-\sqrt 2
\end{array} \right )^{-1}\left (\begin{array}{c} x_n \\
y_n
\end{array} \right ), \ \ \ s_n=\left (\begin{array}{c} 2-\sqrt 2 \\
2+\sqrt 2
\end{array} \right )\frac{c-1}{4n},
$$
and $|r_n|\le \frac{G}{n^2}$ for some $G$. Looking for the
solution in the form
$$
 \varphi _n =A^n c_n
$$
we have the following recurrent formula for $c_n$:
$$
c_{n+1}=A^{-n}(s_n+r_n)+c_n.
$$
Integrating one gets
$$
c_{n+1}=c_1+A^{-2}(s_1+r_1)+A^{-3}(s_2+r_2)+ ... +
A^{-1-n}(s_n+r_n).
$$
For components of $\varphi_n=(a_n,b_n)^T$ this implies
$$
\begin{array}{l}
a_n=\lambda_1^n \left(a_1 + \frac{(2-\sqrt 2)(c-1)}{4}\left(
\frac{1}{\lambda_1^2}+\frac{1}{2\lambda_1^3}+...+\frac{1}{(n-1)\lambda_1^n}
\right)   +\left(
\frac{G_1}{\lambda_1^2}+\frac{G_2}{2^2\lambda_1^3}+...+\frac{G_{n-1}}{(n-1)^2\lambda_1^n}
\right) \right)=
\\
= a_1\lambda_1^n + \frac{(2-\sqrt 2)(c-1)}{4}\left(
\frac{1}{n-1}+\frac{\lambda_1}{(n-2)}+...+\lambda_1^{n-2}
\right)   +\left(
\frac{G_1}{\lambda_1^2}+\frac{G_2}{2^2\lambda_1^3}+...+\frac{G_{n-1}}{(n-1)^2\lambda_1^n}
\right)
\end{array}
$$
with some limited sequence $G_n$: $|G_n|\le G$. The first sum,
corresponding to $s_n$, is estimated as follows:
$$
\begin{array}{l}
\frac{1}{n-1}\left(1+\frac{n-1}{n-2}\lambda_1+\frac{n-1}{n-3}\lambda_1^2...+(n-1)\lambda_1^{n-2}\right)=\\
=\frac{1}{n-1}\left(1+\lambda_1+\lambda_1^2...+\lambda_1^{n-2}\right)+
\frac{\lambda_1}{n-1}\left(\frac{1}{n-2}+\frac{2}{n-3}\lambda_1...+(n-2)\lambda_1^{n-3}\right)=\\
=\frac{1}{n-1}\frac{1-\lambda_1^{n-1}}{1-\lambda_1}+
\frac{\lambda_1}{(n-1)(n-2)}\left(1+\frac{2(n-2)}{n-3}\lambda_1+\frac{3(n-2)}{n-4}\lambda_1^2...+(n-2)^2\lambda_1^{n-3}\right)=\\
=\frac{1}{n-1}\frac{1-\lambda_1^{n-1}}{1-\lambda_1}+
\frac{\lambda_1}{(n-1)(n-2)}F_1(n,\lambda_1)
\end{array}
$$
where
$$
F_1(n,\lambda_1)<(1+2\times 3 \lambda_1+3\times 4
\lambda_1^2+...)=\frac{2}{(1-\lambda_1)^3}-1,
$$
as $|\lambda_1|<1$ and $(n-2)/(n-k)<k$ for $k\le n-1$. The second
sum is estimated by
$$
\begin{array}{l}
\frac{G}{(n-1)^2}\left(1+\frac{(n-1)^2}{(n-2)^2}\lambda_1+\frac{(n-1)^2}{(n-3)^2}\lambda_1^2
+...+\frac{(n-1)^2}{(n-k)^2}\lambda_1^{k-1}+...+(n-1)^2\lambda_1^{n-2}\right)=\\
=\frac{G}{(n-1)^2}(1+\lambda_1+\lambda_1^2+...+\lambda_1^{n-2}+\\
+\frac{(n-1)^2-(n-2)^2}{(n-2)^2}\lambda_1+\frac{(n-1)^2-(n-3)^2}{(n-3)^2}\lambda_1^2
+...+\frac{(n-1)^2-(n-k)^2}{(n-k)^2}\lambda_1^{k-1}+...+((n-1)^2-1)\lambda_1^{n-2})=\\
=\frac{G}{(n-1)^2}\left(\frac{1-\lambda_1^{n-1}}{1-\lambda_1}+
\frac{(2(n-2)+1}{(n-2)^2}\lambda_1
+...+\frac{(2(n-k)(k-1)+(k-1)^2}{(n-k)^2}\lambda_1^{k-1}+...+(2(n-2)+(n-2)^2)\lambda_1^{n-2}\right)\le
\\
\le \frac{G}{(n-1)^2}\left(\frac{1}{1-\lambda_1}+
2\lambda_1(1+2\lambda_1+3\lambda_1^2+...)+\lambda_1(1+2^2\lambda_1+3^2\lambda_1^2+...)
\right)=\frac{G}{(n-1)^2}F_2(\lambda_1).
\end{array}
$$
Summing up we conclude:
\begin{equation}\label{an}
a_n=\frac{(2-\sqrt
2)(c-1)}{4(1-\lambda_1)}\frac{1}{n}+O\left(\frac{1}{n^2}\right).
\end{equation}
For the second component $b_n$ one has
$$
\begin{array}{l}
b_n = b_1\lambda_2^n + \frac{(2+\sqrt 2)(c-1)}{4}\lambda_2^n\left(
\frac{1}{\lambda_2^2}+\frac{1}{2\lambda_2^3}+...+\frac{1}{(n-1)\lambda_2^n}
\right) +\lambda_2^n\left(
\frac{H_1}{\lambda_2^2}+\frac{H_2}{2^2\lambda_2^3}+...+\frac{H_{n-1}}{(n-1)^2\lambda_2^n}
\right)
\end{array}
$$
with some limited sequence $H_n$: $|H_n|\le H$. The first sum in
the previous formula is estimated as
$$
\begin{array}{l}
\lambda_2^{n-1}\int_0^{\frac{1}{\lambda_2}}
(1+x+x^2...+x^{n-2})dx=\lambda_2^{n-1}\int_0^{\frac{1}{\lambda_2}}
\frac{1-x^{n-1}}{1-x}dx=\\
=\lambda_2^{n-1}\left(\int_0^{\frac{1}{\lambda_2}}
\frac{1}{1-x}dx-\int_0^{\frac{1}{\lambda_2}} (x^{n-1}+x^n+...)dx
\right)=\\
=\lambda_2^{n-1}\left(-\ln
(1-x)\mid_0^{\frac{1}{\lambda_2}}-\frac{1}{\lambda_2^n}\left(\frac{1}{n}+\frac{1}{(n+1)\lambda_2}+\frac{1}{((n+2)\lambda_2^2}...
\right) \right)=\\
=\lambda_2^{n-1}\ln\frac{\lambda_2}{\lambda_2-1}-\frac{1}{n\lambda_2}\left(1+\frac{n}{(n+1)\lambda_2}+\frac{n}{((n+2)\lambda_2^2}...
\right)=\\
=\lambda_2^{n-1}\ln\frac{\lambda_2}{\lambda_2-1}-\frac{1}{n\lambda_2}\left(
\left(1+\frac{1}{\lambda_2}+\frac{1}{\lambda_2^2}+... \right)-
\frac{1}{\lambda_2}\left(\frac{1}{(n+1)}+\frac{2}{(n+2)\lambda_2}+...+\frac{k}{(n+k)\lambda_2^{k-1}}+...
\right)\right)=\\
=\lambda_2^{n}\frac{1}{\lambda_2}\ln\frac{\lambda_2}{\lambda_2-1}-\frac{1}{n(\lambda_2-1)}+\frac{S_n(\lambda_2)}{n(n+1)\lambda_2^2},
\end{array}
$$
Where
$$
\begin{array}{l}
S_n(\lambda_2)=\left(1+\frac{2(n+1)}{(n+2)\lambda_2}+\frac{3(n+1)}{(n+3)\lambda_2^2}+...+\frac{k(n+1)}{(n+k)\lambda_2^{k-1}}+...
\right)<\left(1+\frac{2}{\lambda_2}+\frac{3}{\lambda_2^2}+...+\frac{k}{\lambda_2^{k-1}}+...
\right)=\frac{\lambda_2^2}{(\lambda_2-1)^2}.
\end{array}
$$
For the second sum in the formula for $b_n$ one has
$$
\begin{array}{l}
\frac{1}{\lambda_2^2}\left(
H_1+\frac{H_2}{2^2\lambda_2}+...+\frac{H_{n-1}}{(n-1)^2\lambda_2^{n-2}}
\right)=F_3(\lambda_2)-\frac{1}{\lambda_2^2}\left(
\frac{H_n}{n^2\lambda_2^{n-1}}+\frac{H_{n+1}}{(n+1)^2\lambda_2^{n}}+...+\frac{H_{n+k}}{(n+k)^2\lambda_2^{n+k-1}}+...
\right)
\end{array}
$$
where
$$
\begin{array}{l}
F_3(\lambda_2)=\frac{1}{\lambda_2^2}\left(
H_1+\frac{H_2}{2^2\lambda_2}+...+\frac{H_{k-1}}{(k-1)^2\lambda_2^{k-2}}+...
\right)
\end{array}
$$
and the second sum is estimated from above as
$$
\begin{array}{l}
\frac{H}{\lambda_2^{n+1}n^2}\left(1+
\frac{n^2}{(n+1)^2\lambda_2}+\frac{n^2}{(n+2)^2\lambda_2^2}+...+\frac{n^2}{(n+k)^2\lambda_2^k}+...
\right)<\frac{H}{\lambda_2^{n+1}n^2}\left(1+
\frac{1}{\lambda_2}+\frac{1}{\lambda_2^2}+\frac{1}{\lambda_2^2}+...
\right)=\frac{H}{\lambda_2^{n}n^2(\lambda_2-1)}.
\end{array}
$$
Finally
\begin{equation}\label{bn}
b_n=\lambda_2^n\left(b_1+\frac{(2+\sqrt
2)(c-1)}{4\lambda_2}\ln\frac{\lambda_2}{\lambda_2-1}+
F_3(\lambda_2)\right) -\frac{(2+\sqrt
2)(c-1)}{4(\lambda_2-1)}\frac{1}{n}+O\left(\frac{1}{n^2}\right).
\end{equation}
As $b_n\to 0$ with $n\to \infty$ and $\lambda_2>1$ one deduces
that the coefficient by $\lambda_2^n$ vanishes. For original
variables $x_n, y_n$ the found asymptotics (\ref{an}),(\ref{bn})
have especially simple form:
\begin{equation}\label{xnyn}
\left(
\begin{array}{l}
x_n\\
y_n
\end{array}\right)
=\frac{c-1}{2n} \left(
\begin{array}{l}
0\\
1
\end{array}\right)+O\left(\frac{1}{n^2}\right).
\end{equation}
Asymptotic (\ref{xnyn}), the second equation (\ref{t-def}) and
$$
R_{N_0,N_0+n}=R_{N_0,N_0}\prod
\limits_{k=1}^n\frac{2k+(c-1)}{2k-(c-1)}
$$
imply (\ref{as2}). The independence of $K(c)$ on $N_0$ easily
follows from the first equation (\ref{t-def}) and  $x_n\to 0$.
 $\Box$

\section{Discussion and concluding remarks}
Found asymptotics implies
\begin{equation}\label{as4}
\tan\alpha _n \simeq \left(1+\frac{1}{n}\right)^{c-1} \ \ \ {\rm
as} \ \ \ n\to \infty
\end{equation}
for the corresponding solution $u_n=e^{i\alpha_n}$ of
(\ref{dPII}).\\ Further, equation (\ref{c}) allows one to
"integrate" asymptototics (\ref{as2}) to get
\begin{equation}\label{as3}
Z^c(n_0+n,m_0+n) \simeq e^{c \pi i/4}K(c)n^{c} \ \ \ {\rm as} \ \
\ n\to \infty.
\end{equation}
Thus the circles of $Z^c$ not only cover the whole infinite
sector with the angle $c\pi/2$ but the circle centers and
intersection points mimic smooth map $z\to z^c$ also
asymptotically. Morover, $R(z)$, being analogous to $|f'(z)|$ of
the corresponding smooth map, has the "right" asymptotics as well.

Discrete map $Z^c$ is defined via constraint (\ref{c}) which is
isomonodromy condition for some linear equation. This seems to be
rather far-fetched approach. It would be more naturally to define
$Z^c$ via circle patterns in a pure geometrical way: square grid
$Z^c$ with $0 <c<2$ is an embedded infinite square grid circle
pattern, the circles $C(z)$ being labeled by
$$
{\mathbb V}=\{z=N+iM:\ N,M \in {\bf Z^2}, M \ge |N| \},
$$
satisfying the following conditions:\\
1) the circles cover the infinite sector with the angle
$c\pi/2$,\\
2) the centers of the border circles $C(N+iN)$ and $C(-N+iN)$ lie
on the borders of this sector.
\begin{conjecture}\label{geodefcon}
Up to re-scaling there is unique square grid $Z^c$.
\end{conjecture}

\begin{figure}[th]
\begin{center}
\epsfig{file=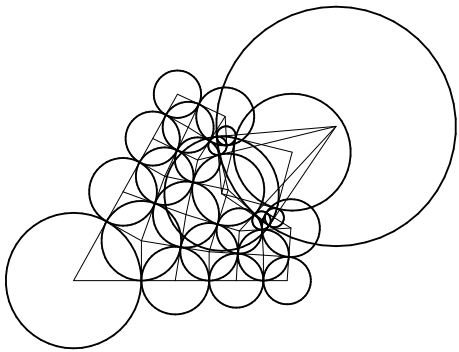,width=60mm} \caption{Non-immersed
discrete conformal map} \label{correct-wrong}
\end{center}
\end{figure}

A naive method to look for so defined discrete $Z^c$ is to start
with some equidistant  $f_{n,0}\in \mathbb R$, $f_{0,m}\in i^{c
}\mathbb R$:
$$
|f_{2n,0}-f_{2n+1,0}|=|f_{2n,0}-f_{2n-1,0}|=|f_{0,2n}-f_{0,2n+1}|=|f_{0,2n}-f_{0,2n-1}|
$$
and then compute $f_{n,m}$ for any $n,m > 0$ using equation
(\ref{q}). Such $f_{n,m}$ determines some circle patterns as
(\ref{equi}) holds true.    But so determined map has a not very
nice behavior (see an example in Fig.~\ref{correct-wrong}).
For $c=2/k,\ k\in \mathbb N$ the proof of Conjecture
\ref{geodefcon} easily follows from the rigidity results obtained
in \cite{H}.

Infinite embedded circle patterns define some (infinite) convex
ideal polyhedron in $\mathbb H^3$ and the group generated by
inversions in its faces.  The known results (see for example
\cite{S}) imply the conjecture claim for rational $c$, but seem
to be inapplicable for irrational as the group generated but the
circles and the sector borders is not discrete any more.
Unfortunately, the rigidity results for finite polyhedra \cite{Ri}
does not seem to be carried over to infinite case by induction.

One can also consider solutions of (\ref{q}) subjected to
(\ref{c}) where $n$ and $m$ are not integer. It turned out that
there exist initial data, so that the corresponding solution
define immersed circle patterns. (The equations for radii are the
same and therefore compatible. Existence of positive solution is
provable by the arguments for the corresponding discrete
Painlev\'e equation as in \cite{AB}.) It is natural to call such
circle patterns discrete map  $z\to (z+z_0)^c$. In this case it
can not be defined pure geometrically as the centers of  border
circles are not collinear.

One can discard the restriction $0 <c<2$  (as well as some
circles $C(z)$ with small $z$) to define discrete $Z^K$ for
natural $K>2$ (see \cite{AB}). All the asymptotic results
obtained so far can be carried over on this case as well as on
hexagonal $Z^c$ with constant intersection angles defined in
\cite{BH} since the governing equations are essentially the same
(\cite{ABh}).

\section{Acknowledgements}

The author thanks B.Apanasov, A.Bobenko and Yu.Suris for useful
discussions.

\end{document}